\documentclass[a4paper,12pt,pagesize]{scrartcl}
\pdfoutput=1

\usepackage[utf8]{inputenc}
\usepackage[english,ngerman]{babel}
\usepackage[T1]{fontenc}
\usepackage{lmodern}
\usepackage{microtype}
\usepackage{amsmath}
\usepackage{amssymb}
\usepackage{hyperref}
\usepackage{hypcap}
\usepackage[all]{xy}
\usepackage{tikz}
\usepackage[hyperref,amsmath,thmmarks,thref]{ntheorem}
\usepackage{dsfont}

\title{Homological stability for classical groups}
\author{Jan Essert}
\date{\today}

\xdefinecolor{myblue}{rgb}{0.1,0.1,0.6}
\xdefinecolor{mygreen}{rgb}{0.1,0.4,0.3}

\hypersetup{
pdftitle={Homological stability for classical groups},
pdfauthor={Jan Essert},
pdfkeywords={Homological stability, Classical groups, Buildings, Group homology, Group theory, K-theory, Opposition complexes},
colorlinks=true,
linkcolor=myblue,
citecolor=mygreen,
urlcolor=myblue
}

\newcommand{\Z}{\mathds Z}

\newcommand{\one}{\mathds 1}

\newcommand{\cA}{\mathcal A}

\newcommand{\cH}{\mathcal H}

\DeclareMathOperator{\id}{id}

\DeclareMathOperator{\Gl}{GL}
\DeclareMathOperator{\Sl}{SL}
\DeclareMathOperator{\SO}{SO}
\DeclareMathOperator{\Orth}{O}
\DeclareMathOperator{\U}{U}

\DeclareMathOperator{\Sp}{Sp}

\DeclareMathOperator{\rank}{rank}
\DeclareMathOperator{\type}{type}
\DeclareMathOperator{\lk}{lk}
\DeclareMathOperator{\st}{st}
\DeclareMathOperator{\Cone}{Cone}
\DeclareMathOperator{\Tot}{Tot}

\newcommand{\coloneq}{\mathrel{\mathop :}=}

\theoremstyle{plain}
\newtheorem{Definition}{Definition}[section]
\newtheorem{Theorem}[Definition]{Theorem}
\newtheorem{Lemma}[Definition]{Lemma}
\newtheorem{Proposition}[Definition]{Proposition}

\newtheorem*{RelSpecTheorem}{Theorem \protect\ref{th:stability_pair_spectral_sequence}}
\newtheorem*{SlnTheorem}{Theorem \protect\ref{th:sln_result}}
\newtheorem*{UnTheorem}{Theorem \protect\ref{th:un_result}}
\newtheorem*{SOnTheorem}{Theorem \protect\ref{th:son_result}}

\theoremclass{Theorem}
\theoremstyle{break}

\theoremstyle{nonumberplain}
\theorembodyfont{\rm}
\newtheorem{Remark}{Remark}

\newtheorem{Example}{Example}

\newtheorem{Construction}{Construction}
\theoremsymbol{\ensuremath\Box}
\newtheorem{Proof}{Proof}

\addtokomafont{paragraph}{\rm\bf}
\addtokomafont{subparagraph}{\rm\it}

\begin{document}
\selectlanguage{english}
\maketitle
\begin{abstract}
  Associated to every group with a weak spherical Tits system of rank $n+1$ with an appropriate rank $n$ subgroup, we construct a relative spectral sequence involving group homology of Levi subgroups of both groups. Using the fact that such Levi subgroups frequently split as semidirect products of smaller groups, we prove homological stability results for unitary groups over division rings with infinite centre as well as for special linear and special orthogonal groups over infinite fields.
\end{abstract}

\addtocounter{section}{-1}
\section{Introduction}

Homological stability is the following question: Given an infinite series of groups $G_n$, such as the general linear groups $\Gl_n$, we consider the sequence of inclusions
\[
	G_1 \hookrightarrow G_2 \hookrightarrow G_3 \hookrightarrow \cdots.
\]
Then, if we apply group homology of a fixed degree, does the corresponding sequence of homology modules stabilise eventually?

This is an old question and there are many interesting results for various series of classical groups, usually over rings of finite stable rank. An overview of results in this area can be found in \cite[Chapter 2]{Knu:HLG:01} and we will also provide references to the best known results for specific series of groups. Although the method of proof is usually based on a common idea, the action of the larger group on a highly connected simplicial complex, all proofs known to the author are tailored to specific series of groups.

In this paper, we present a general method to prove homological stability, valid for all groups with weak spherical Tits systems, that is, groups acting strongly transitively on possibly weak spherical buildings. We then use this method to prove homological stability for various series of classical groups over division rings, usually improving the stability range previously known from work with larger classes of rings.

The method is based on the observation that the simplicial complexes used by Charney in \cite{Cha:HSD:80} and \cite{Cha:gtV:87} and by Vogtmann in \cite{Vog:HSO:79} and \cite{Vog:SPH:81} are closely related to the theory of buildings --- they are the \emph{opposition complexes} studied by von Heydebreck in \cite{vH:HPC:03}. The opposition complexes admit Levi subgroups as vertex stabilisers. Using these complexes, we construct a spectral sequence involving relative group homology of Levi subgroups.

For the groups we consider, the Levi subgroups split as direct or semidirect products of smaller groups, both of the series of groups we consider and, interestingly, of \emph{general linear groups}. Using strong stability results for general linear groups, we can hence prove stability results for various series of groups.

\paragraph{Homological stability results} The method outlined above has originally been used by Charney in \cite{Cha:HSD:80} to prove homological stability of general and special linear groups, but yielding a comparatively weak stability range. For special linear groups, however, it is an interesting observation that terms involving general linear groups appear in the spectral sequence. This allows us to apply a strong theorem by Sah in \cite{Sah:HcL:86} on homological stability for general linear groups to prove homological stability for special linear groups.

\begin{SlnTheorem}[Homological stability of special linear groups]
  If $D$ is an infinite field, then $n\geq 2k-1$ implies
  \[
	  H_k(\Sl_{n+1}(D),\Sl_{n}(D);\Z)=0.
  \]
\end{SlnTheorem}

\noindent For fields of characteristic zero, there is a far better result by Hutchinson and Tao in \cite{HT:HSS:08} with stability range $n\geq k$. Before our work, the best result known to the author applicable to other infinite fields was a result by van der Kallen in \cite{vdK:HSL:80} for rings with stable rank $1$. It guarantees a stability range of $n\geq 2k$.

Vogtmann originally used a version of the construction in this paper to prove homological stability for orthogonal and symplectic groups in \cite{Vog:HSO:79} and \cite{Vog:SPH:81}. Here, we investigate the general situation of unitary groups associated to a hermitian form of Witt index $n+1$ on a vector space $V$. This vector space then splits non-canonically as an orthogonal sum of a hyperbolic module $\cH_{n+1}$ and an anisotropic complement $W$. We consider the unitary group induced on the subspace $\cH_n\perp W$ and ask for homological stability. Again, the spectral sequence we consider has terms involving the relative homology of general linear groups. We can hence apply Sah's theorem again to obtain

\begin{UnTheorem}[Homological stability of unitary groups]
  For a division ring $D$ with infinite centre, the relative homology modules
  \[
    H_k\bigl(\U(\cH_{n+1}\perp W), \U(\cH_n\perp W);\Z\bigr)
  \]
  vanish for $n\geq 2$ if $k=1$ and for $n\geq k\geq 2$. If the centre of $D$ is finite, relative homology vanishes for $n\geq 2k$.
\end{UnTheorem}

\noindent This is an improvement over the results by Mirzaii and van der Kallen in \cite{MaB:HSU:02} and \cite{Mir:HSU:05}, where homological stability for unitary groups with stability range $n\geq k+1$ has already been proved. Their result is valid for a much larger class of rings, namely local rings with infinite residue field, but only for the case of maximal Witt index, that is for $W=\{0\}$.

The following strong result will also be proved using this method.

\begin{SOnTheorem}[Homological stability of special orthogonal groups]
  For an infinite field $D$, we have
  \[
    H_k\bigl(\SO_{n+1,n+1}(D),\SO_{n,n}(D);\Z\bigr) =0
  \]
  for $n\geq 2$ if $k=1$ and for $n\geq k\geq 2$. If $D$ is a finite field, then the relative homology groups vanish for $n\geq 2k$.
\end{SOnTheorem}

\paragraph{The construction of a relative spectral sequence} We give an outline of the method used to prove these results and we state the main theorem. Consider a group $G$ with a weak spherical Tits system of rank \mbox{$n+1$}, contained in an infinite series of groups for which we prove homological stability. We enumerate the type set $I=\{i_1,\ldots,i_{n+1}\}$ arbitrarily. For $1\leq p\leq n+1$, denote by $L_p$ certain Levi subgroups of $G$ of type $I\backslash\{i_p\}$.

For the applications discussed in the previous section, the group $G$ is of type $A_{n+1}$ or $C_{n+1}$ with a linear ordering of the type set. The resulting Coxeter diagrams of $G$ and $L_p$ are illustrated in the following picture.
\begin{center}
\begin{tikzpicture}[font=\small]
  \node (G) at (-1,.8) {$G$};
  \node (L) at (-1,0) {$L_p$};
  \foreach \y in {0,.8} {
  \foreach \x in {0,1,2,3,4,6,7,8,9,10} { \fill (\x,\y) circle (.7mm);}
  \draw (1,\y) -- (2,\y);
  \draw[dotted] (2,\y) -- (3,\y);
  \draw (3,\y) -- (4,\y);
  \draw (6,\y) -- (7,\y);
  \draw[dotted] (7,\y) -- (8,\y);
  \draw (8,\y) -- (10,\y);
  \draw (0,\y + .05) -- (1,\y + .05);
  \draw[dashed]  (0,\y - .05) -- (1,\y -.05);
  }
  \draw (4,.8) -- (6,.8);
  \fill (5,.8) circle (.7mm);
  \node (0) at (0,-.6) {$1$}; \node (1) at (1,-.6) {$2$}; \node (3) at (4,-.6) {$p-1$}; \node (4) at (5,-.63) {$p$}; \node (5) at (6,-.6) {$p+1$}; \node (7) at (9,-.62) {$n$}; \node (8) at (10,-.6) {$n+1$};
\end{tikzpicture}
\end{center}
We choose a subgroup $G'\leq L_{n+1}$ of type $I\backslash \{i_{n+1}\}$ and write $L'_p=L_p\cap G'$. Again, in the concrete applications, we have the following situation.
\begin{center}
\begin{tikzpicture}[font=\small]
  \node (G) at (-1,.8) {$G'$};
  \node (L) at (-1,0) {$L'_p$};
  \foreach \y in {0,.8} {
  \foreach \x in {0,1,2,3,4,6,7,8,9} { \fill (\x,\y) circle (.7mm);}
  \draw (1,\y) -- (2,\y);
  \draw[dotted] (2,\y) -- (3,\y);
  \draw (3,\y) -- (4,\y);
  \draw (6,\y) -- (7,\y);
  \draw[dotted] (7,\y) -- (8,\y);
  \draw (8,\y) -- (9,\y);
  \draw (0,\y + .05) -- (1,\y + .05);
  \draw[dashed]  (0,\y - .05) -- (1,\y -.05);
  }
  \draw (4,.8) -- (6,.8);
  \fill (5,.8) circle (.7mm);
  \node (0) at (0,-.6) {$1$}; \node (1) at (1,-.6) {$2$}; \node (3) at (4,-.6) {$p-1$}; \node (4) at (5,-.63) {$p$}; \node (5) at (6,-.6) {$p+1$}; \node (7) at (9,-.6) {$n$};
  \node (empty) at (10,0) {\phantom{$n+1$}};
\end{tikzpicture}
\end{center}
Using a filtration of the opposition complex by types of vertices, we construct two exact chain complexes of $G$- and $G'$-modules. From these chain complexes, we obtain a spectral sequence involving relative homology of Levi subgroups with coefficient modules $M_p$, which are top-dimensional homology modules of opposition complexes of type $\{i_1,\ldots,i_{p-1}\}$, except for $M_1=\Z$.

\begin{RelSpecTheorem}[Relative spectral sequence]
  There is a spectral sequence with first page
  \[
  E^1_{p,q}=\begin{cases}
    H_q(G,G';\Z) & p=0 \\
    H_q(L_p,L'_p;M_p) & 1\leq p \leq n\\
    H_q(L_{n+1},G';M_{n+1}) & p=n+1
  \end{cases}
  \]
  which converges to zero.
\end{RelSpecTheorem}

This can be used to prove homological stability for groups of type $A_{n+1}$ and $C_{n+1}$ in the following way: We want to prove that $H_q(G,G';\Z)$ vanishes for all $q$ smaller than a given $k$. Hence we must show that $H_q(L_p,L'_p;M_p)$ vanishes for $p+q\leq k+1$. For $2\leq p\leq n-1$ the Levi subgroups, having disconnected diagrams, usually split as direct or semidirect products of two groups whose types belong to the connected components of the diagrams.
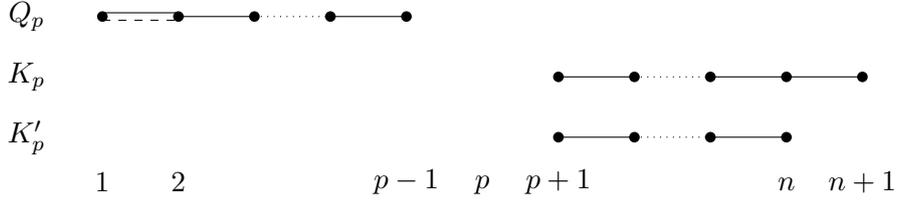
\begin{center}
\begin{tikzpicture}[font=\small]
  \node (Qp) at (-1,1.6) {$Q_p$};
  \node (K) at (-1,0.8) {$K_p$};
  \node (Kp) at (-1,0) {$K'_p$};
  \foreach \x in {0,1,2,3,4} { \fill (\x,1.6) circle (.7mm);}
  \draw (1,1.6) -- (2,1.6);
  \draw[dotted] (2,1.6) -- (3,1.6);
  \draw (3,1.6) -- (4,1.6);
  \draw (0,1.6 + .05) -- (1,1.6 + .05);
  \draw[dashed]  (0,1.6 - .05) -- (1,1.6 -.05);
  \foreach \y in {0,.8} {
  \foreach \x in {6,7,8,9} { \fill (\x,\y) circle (.7mm);}
  \draw (6,\y) -- (7,\y);
  \draw[dotted] (7,\y) -- (8,\y);
  \draw (8,\y) -- (9,\y);
  }
  \draw (9,.8) -- (10,.8);
  \fill (10,.8) circle (.7mm);
  \node (0) at (0,-.6) {$1$}; \node (1) at (1,-.6) {$2$}; \node (3) at (4,-.6) {$p-1$}; \node (4) at (5,-.63) {$p$}; \node (5) at (6,-.6) {$p+1$}; \node (7) at (9,-.62) {$n$}; \node (8) at (10,-.6) {$n+1$};
\end{tikzpicture}
\end{center}
This means that there are groups $Q_p$, $K_p$, $K'_p$ of types $\{i_1,\ldots,i_{p-1}\}$, $\{i_{p+1},\ldots,i_{n+1}\}$ and $\{i_{p+1},\ldots,i_n\}$, respectively, such that $L_p=Q_p\ltimes K_p$ and $L'_p=Q_p\ltimes K'_p$. The modules $M_p$ are constructed in such a fashion that the groups $K_p$ and $K'_p$ act trivially on $M_p$. If we know that relative \emph{integral} homology of the subgroups $K_p$ and $K'_p$  vanishes, we can produce zeroes in the spectral sequence $E^1_{p,q}$ by using a relative Lyndon\slash Hochschild-Serre spectral sequence.

The structure of the Levi subgroups and the corresponding semidirect product decompositions depend on the specific series of groups, but note that we always need relative integral homology of groups of type $A_*$! For special linear groups over fields and for unitary groups over division rings, these subgroups of type $A_*$ are general linear groups. Hence, as mentioned above, we can use strong results on the homological stability of general linear groups to obtain homological stability of these series of groups.

Variations of this method using an appropriate type filtration can probably be used to prove homological stability results for different series of reductive groups. As mentioned above, this could be used to study groups of type $E_6$, $E_7$ or $E_8$, or compare group homology of groups of different types. 

\paragraph{Possible applications} The relative spectral sequence could be useful for various applications the author did not have time to look at. One example could be low-dimensional homological stability for groups of types $E_6$, $E_7$ and $E_8$. Additionally, one could try to compare group homology of groups of different types using this method. In particular, relations between algebraic K-theory and hermitian K-theory could also be studied by choosing a different type enumeration of a group of type $C_{n+1}$, forcing $G'$ to be of type $A_n$ instead of type $C_n$. Finally, homological stability results for all reductive algebraic groups should be possible, albeit with a rather weak stability range. 

\paragraph{} This paper is structured as follows. In the first part, we introduce the required concepts very briefly and give a general procedure to construct relative spectral sequences. At the end of this section, we apply this procedure to obtain a relative Lyndon\slash Hochschild-Serre spectral sequence, which is required to decompose homology of Levi subgroups later on.

In the second part, we introduce the groups we work with and their associated opposition complexes. We construct an exact chain complex, analogously to the construction of cellular homology, and use it to prove Theorem \ref{th:stability_pair_spectral_sequence}.

Finally, in the third part, we consider explicit series of groups, decompose the Levi subgroups appropriately and apply Theorem \ref{th:stability_pair_spectral_sequence} to prove homological stability inductively.

\paragraph{} The author would like to thank Linus Kramer for many discussions, suggestions and encouragement, as well as Ruth Charney and Karen Vogtmann for their help. We also thank Stefan Witzel for very useful discussions. Proposition \ref{prop:filtrated_opposition_complexes_are_spherical} is largely due to him.
The author was supported by the \selectlanguage{ngerman}\emph{Graduiertenkolleg: \glqq Analytische Topologie und Me\-ta\-geo\-me\-trie\grqq}\selectlanguage{english} while working on this topic. This work is part of the author's doctoral thesis \cite{Ess:BGL:10} at the Universität Münster.

\section{Homology}

The aim of this part is to give a general construction for relative spectral sequences. Using this, we will construct a relative Lyndon\slash Hochschild-Serre spectral sequence. The relative construction will also be applied in the second part to prove the main theorem. We will assume the reader to be familiar with the concept of spectral sequences. For a textbook on this topic, see \cite{McC:UGS:01}. A very good brief introduction to the subject is \cite{Cho:Ych:06}.

\subsection{Group homology}

The reader is also assumed to be familiar with group homology. The standard reference is of course \cite{Bro:CoG:82}. Results from this book will be referred to throughout this paper. In this section, we describe all non-standard terminology we will use.

Throughout this section, let $G$ be any group, we write $\Z G$ for the group ring over $G$, and we define a \emph{$G$-module} to be a left module over $\Z G$. The standard resolution of $\Z$ over $\Z G$, as defined in \cite[I.5]{Bro:CoG:82} will be denoted by $F_*(G)$. The group homology of $G$ is defined to be the homology of this chain complex. Since \cite{Bro:CoG:82} does not contain relative homology, we define it here.

\begin{Definition}
  Consider a group $G$ and a subgroup $G'\leq G$. Let $M$ be a $G$-module. The \emph{relative group homology} $H_*(G,G';M)$ is defined to be the homology of the quotient complex $F_*(G)\otimes_G M / F_*(G') \otimes_{G'} M$.
\end{Definition}

\noindent We will denote the \emph{mapping cone chain complex} of a chain map $\varphi_*:C'_*\rightarrow C_*$ by $\Cone_*(\varphi)$, for details see \cite[1.5]{Wei:IHA:94}.

Extensive use will be made of double complexes and their homology, as in \cite[VII.3]{Bro:CoG:82}. We will denote the \emph{total complex} associated to a double complex $D_{p,q}$ by $\Tot(D)$.

In addition to the tensor complex $F\otimes_G C$ of two chain complexes of $G$-modules, we will also use the \emph{transposed tensor complex} denoted by
\[
(F\otimes^T_G C)_k \coloneq \bigoplus_{p+q=k} F_p \otimes_G C_q
\]
with differential induced by
\[
\partial^{F\otimes^T_G C}(f_p \otimes_G c_q) = (-1)^q \partial^F(f_p) \otimes_G c_q + f_p \otimes_G \partial^C(c_q).
\]
This is obviously the total complex of the double complex $(F_q \otimes_G C_p)_{p,q}$. The following Lemma is then a simple calculation.

\begin{Lemma}\label{lem:transposed_double_complex}
    The chain map $F\otimes_G C\rightarrow F\otimes^T_G C$ given on the basis by
    \[
    f_p \otimes_G c_q \mapsto (-1)^{pq} f_p \otimes_G c_q
    \]
    induces isomorphisms on homology $H_*(F\otimes_G C) \cong H_*(F\otimes^T_G C)$.
\end{Lemma}

\subsection{Relative spectral sequences}

In this section, we will discuss a procedure to obtain relative spectral sequences from both the spectral sequence associated to a double complex as well as from the Lyndon\slash Hochschild-Serre spectral sequence. Consider the following situation:

Fix a group $G$ and a subgroup $G'$. Let $F'$, $C'$ and $F$, $C$ be chain complexes of $G'$- and $G$-modules, respectively. Consider the two tensor product double complexes $(F_p \otimes_G C_q)_{p,q}$ and $(F'_p \otimes_{G'} C'_q)_{p,q}$. Assume that there is map of double complexes
\[
i: F'\otimes_{G'} C' \rightarrow F\otimes_G C.
\]
We denote the induced maps on the vertical and horizontal chain complexes by
    \begin{align*}
        i_{p,\bullet}: F'_p \otimes_{G'} C'_\bullet & \rightarrow F_p \otimes_G C_\bullet \\
        i_{\bullet,q}: F'_\bullet \otimes_{G'} C'_q & \rightarrow F_\bullet \otimes_G C_q
    \end{align*}
    For every $q$, we then consider the mapping cone chain complexes $\Cone_*(i_{\bullet,q})$ and $\Cone_*(i_{q,\bullet})$. It is a simple calculation to see that
    \[
    D_{p,q}=\Cone_p(i_{\bullet,q})\qquad\text{and}\qquad D^T_{p,q}=\Cone_p(i_{q,\bullet})
    \]
    are actually double complexes with respect to the cone differentials and the differentials
	\begin{align*}
        \partial: (F'_{p-1} \otimes_{G'} C'_q) \oplus (F_p \otimes_G C_q) &\longrightarrow (F'_{p-1} \otimes_{G'} C'_{q-1}) \oplus (F_p \otimes_G C_{q-1}) \\
        f' \otimes_{G'} c' + f \otimes_G c &\longmapsto f' \otimes_{G'} \partial^{C'}(c') + f\otimes_G\partial^C(c) \\
        \partial^T: (F'_q \otimes_{G'} C'_{p-1}) \oplus (F_q \otimes_G C_p) &\longrightarrow (F'_{q-1} \otimes_{G'} C'_{p-1}) \oplus (F_{q-1} \otimes_G C_p) \\
        f' \otimes_{G'} c' + f \otimes_G c &\longmapsto \partial^{F'}(f') \otimes_{G'} c' + \partial^F (f)\otimes_G c
	\end{align*}
	induced by the differentials $\partial^{C'}$ and $\partial^C$, respectively $\partial^{F'}$ and $\partial^F$.

\begin{Theorem}\label{th:relative_spectral_sequence}
  There are two spectral sequences corresponding to $D$ and $D^T$ satisfying
	\[
    \begin{array}{c} E^1_{p,q} = H_q\bigl(\Cone_*(i_{\bullet,p})\bigr) \\
        L^1_{p,q} = H_q\bigl(\Cone_*(i_{p,\bullet})\bigr) \end{array} \Rightarrow H_{p+q}\bigl(\Cone_*(i)\bigr).
	\]
    The differentials on the first pages are induced by $\pm\partial$ and $\pm\partial^T$, respectively.
\end{Theorem}
\begin{Proof}
    The first spectral sequence is constructed as in \cite[VII.3]{Bro:CoG:82} using the double complex $D$. We know that
    \[
        E^1_{p,q} = H_q\bigl(\Cone_*(i_{\bullet,p})\bigr) \Rightarrow H_{p+q}(\Tot(D)).
    \]
    But the formation of the mapping cone and of the total complex commutes in the following sense: On the one hand
    \[
    \Tot(D)_k = \bigoplus_{p+q=k} (F'_{p-1} \otimes_{G'} C'_q) \oplus (F_p \otimes_G C_q).
    \]
	On the other hand
	\[
    \Cone_k(i:(F'_*\otimes_{G'} C'_*)\rightarrow (F_* \otimes_G C_*)) = \bigoplus_{p+q=k-1}(F'_p \otimes_{G'} C'_q) \oplus \bigoplus_{p+q=k}(F_p \otimes_G C_q)
	\]
	Hence the modules of the two chain complexes coincide. For the boundary maps, we have
	\begin{align*}
        \partial^{\Cone(i)}\Bigl(\sum_{p+q=k-1}f'_p\otimes_{G'} c'_q + \sum_{p+q=k}f_p\otimes_G c_q \Bigr) =
        &- \sum_{p+q=k-1} \bigl(\partial f'_p \otimes_{G'} c'_q + (-1)^p f'_p\otimes_{G'}\partial c'_q \bigr) +\\
        &+ i \bigl( \sum_{p+q=k-1} f'_p \otimes_{G'} c'_q \bigr) +\\
		&+ \sum_{p+q=k} \bigl(\partial f_p\otimes_G c_q + (-1)^p f_p\otimes_G\partial c_q\bigl).
    \end{align*}
		and
    \begin{multline*}
        \partial^{\Tot(D)}\Bigl(\sum_{p+q=k} (f'_{p-1}\otimes_{G'} c'_q + f_p\otimes_G c_q)\Bigr) =
        \sum_{p+q=k} \bigl( -\partial f'_{p-1}\otimes_{G'} c'_q + i(f'_{p-1}\otimes_{G'} c'_q) + \partial f_p\otimes_G c_q + \\
        \qquad\quad+ (-1)^p f'_{p-1}\otimes_{G'} \partial c'_q + (-1)^p f_p\otimes_G \partial c_q\bigr).
	\end{multline*}
	which are also easily seen to be identical.

    For the second spectral sequence, we use the procedure from \cite[VII.3]{Bro:CoG:82} again, this time using the double complex $D^T$. We denote it by $L$ to make the notation consistent in the following. Note that we can also write $D^T_{p,q}\cong\Cone_p(i_{\bullet,q}^T)$, where $i^T: (F'_q\otimes_{G'}^T C'_p)_{p,q} \rightarrow (F_q\otimes_G^T C_p)_{p,q}$ is the map induced on the transposed double complexes. By applying the above calculation to $i^T$, we obtain
	\[
	L^1_{p,q} = H_q\bigl(\Cone_*(i_{p,\bullet})\bigr) \Rightarrow H_{p+q}(\Tot(D^T)) \cong H_{p+q}(\Cone_*(i^T)),
	\]
    We have seen in Lemma \ref{lem:transposed_double_complex} that there is a chain map inducing isomorphisms
	\[
	H_*(F \otimes_G C) \cong H_*(F \otimes_G^T C)
	\]
    on homology. Using these maps, one can easily construct a chain map $\Cone(i)\rightarrow \Cone(i^T)$ inducing an isomorphism on homology by the long exact sequence associated to the mapping cone and the 5-Lemma.
\end{Proof}

\begin{Remark}
  If $i$ is injective, we have
  \[
  H_{p+q}(\Cone_*(i))\cong H_{p+q}( F\otimes_G C, i(F'\otimes_{G'} C'))
  \]
  by \cite[1.5.8]{Wei:IHA:94}.
\end{Remark}

\noindent Applied to the Lyndon\slash Hochschild-Serre spectral sequence as in \cite[VII.6]{Bro:CoG:82}, we obtain

\begin{Theorem}[Relative Lyndon\slash Hochschild-Serre]\label{th:relative_lyndon_hochschild_serre}
	Fix a group $G$ with a normal subgroup $H$, a subgroup $G'\leq G$ and a $G$-module $M$. Write $H'=H\cap G'$. If the quotient $Q=G/H$ is isomorphic to $G'/H'$, there is a spectral sequence
	\[
	L^2_{p,q} = H_p\bigl(Q;H_q(H,H';M)\bigr) \Rightarrow H_{p+q}(G,G';M).
	\]
	If the subgroup $H$ acts trivially on $M$ and $M$ is $\Z$-free, we obtain
	\[
	L^2_{p,q} = H_p\bigl(Q;H_q(H,H';\Z)\otimes_\Z M\bigr) \Rightarrow H_{p+q}(G,G';M).
	\]
\end{Theorem}

\begin{Proof}
	Let $F_*(G)$, $F_*(G')$ and $F_*(Q)$ be the standard resolutions of $\Z$ over $\Z G$, $\Z G'$ and $\Z Q$, respectively. As in the construction of the Lyndon\slash Hochschild-Serre spectral sequence, we consider the double complexes $F_*(Q) \otimes_Q (F_*(G') \otimes_{H'} M)$ and $F_*(Q) \otimes_Q (F_*(G) \otimes_H M)$. Let
    \[
        i : F_*(Q) \otimes_Q ( F_*(G') \otimes_{H'} M) \rightarrow F_*(Q) \otimes_Q (F_*(G) \otimes_H M)
    \]
    be induced by the inclusions. Then apply Theorem \ref{th:relative_spectral_sequence} to obtain a spectral sequence with first page term
	\[
    L^1_{p,q}=H_q(\Cone_*(i_{p,\bullet}))\cong H_q\biggl(\frac{F_p(Q) \otimes_Q (F_*(G) \otimes_{H} M)}{F_p(Q) \otimes_Q (F_*(G') \otimes_{H'} M)}\biggr).
	\]
	The module $F_p(Q)$ is $\Z Q$-free and hence $\Z Q$-flat. We obtain
	\[
	L^1_{p,q}\cong F_p(Q) \otimes_Q H_q\biggl(\frac{F_*(G)\otimes_{H} M}{F_*(G') \otimes_{H'} M}\biggr) \cong F_p(Q) \otimes_Q H_q\biggl(\frac{F_*(H)\otimes_{H} M}{F_*(H') \otimes_{H'} M}\biggr),
	\]
    where the last isomorphism can be seen easily from the associated long exact sequences using the five-lemma. This yields 
    \[
    L^2_{p,q}\cong H_p(Q,H_q(H,H';M)) \Rightarrow H_{p+q}(\Cone_*(i)).
    \]

    \noindent On the other hand, consider the spectral sequence $E$ from Theorem \ref{th:relative_spectral_sequence} and apply \cite[1.5.8]{Wei:IHA:94} to obtain the following description of the first page.
    \[
    E^1_{p,q} \cong H_q \biggl( \frac{F_*(Q) \otimes_Q ( F_p(G) \otimes_H M)}{F_*(Q) \otimes_Q (F_p(G') \otimes_{H'} M)} \biggr) \Rightarrow H_{p+q}(\Cone_*(i)).
    \]
    In the proof of the regular Lyndon\slash Hochschild-Serre spectral sequence as in \cite[VII.6]{Bro:CoG:82}, we have seen that $H_q(F_*(Q) \otimes_Q ( F_p(G) \otimes_H M))=0$ for $q\neq 0$ and that 
    \[
    H_0(F_*(Q) \otimes_Q ( F_p(G) \otimes_H M)) \cong F_p(G)\otimes_G M.
    \]
    The same is true if we replace $G$ and $H$ by $G'$ and $H'$. The map $i$ induces the map
    \[
    H_0(F_*(Q) \otimes_Q ( F_p(G') \otimes_{H'} M)) \rightarrow H_0(F_*(Q) \otimes_Q ( F_p(G) \otimes_H M))
    \]
    which under the above isomorphisms is just the inclusion $F_p(G')\otimes_{G'}M \rightarrow F_p(G)\otimes_G M$. In particular, it is injective. By the long exact sequence for relative homology, we obtain $E^1_{p,q}=0$ for $q\neq 0$ and
   \[
    E^1_{p,0}\cong \frac{F_p(G) \otimes_G M}{F_p(G')\otimes_{G'} M},
    \]
    so the spectral sequence $E$ collapses on the second page and converges to $H_{p+q}(G,G';M)$, which proves the first statement.
    
    For the second statement, since $M$ is trivial, we have $F_k(H) \otimes_H M \cong F_k(H) \otimes_H \Z \otimes_\Z M$, likewise for $H'$. Since $M$ is $\Z$-free, it is in particular $\Z$-flat and the functor $(-\otimes_\Z M)$ is exact, which proves the result.
\end{Proof}

\section{Geometry}

Homological stability proofs usually consider the action of a group on some simplicial complex and then exhibit smaller groups of the same series of groups as stabiliser subgroups. In this part, we will introduce the opposition complex associated to a group with a weak spherical Tits system and construct a filtration of this complex which leads to a relative spectral sequence involving the group and its Levi subgroups. This will be used in the last part to prove homological stability.

\subsection{Spherical buildings}\label{subsec:sph_buildings}

We will briefly recall the basic definitions for Coxeter complexes and spherical buildings. The books by Abramenko and Brown \cite{AB:B:08} and by Ronan \cite{Ron:LoB:89} are excellent references, where all of the material of this section can be found.

At the end of this section, we will illustrate all of these definitions in the concrete case of the projective space over a division ring.

\begin{Definition}
    Let $I$ be a finite set. A \emph{Coxeter matrix} $M=(m_{i,j})_{i,j\in I}$ is a symmetric matrix with entries 1 on the diagonal and with entries in $\{2,3,\ldots,\infty\}$ else.

    We associate to $M$ the \emph{Coxeter diagram}, an edge-labelled graph with vertex set $I$ and with edges between $i$ and $j$ if $m_{i,j}>2$. We label all edges by the corresponding matrix elements $m_{i,j}$. It is customary to omit the label if $m_{i,j}=3,4$ and to draw a double edge for $m_{i,j}=4$.
\end{Definition}

\noindent In this paper, we will need the following two diagrams, where the type set $I=\{i_1,\ldots,i_n\}$ is enumerated linearly as follows
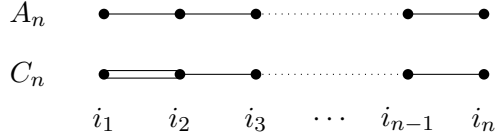
\begin{center}
\begin{tikzpicture}[font=\small]
    \node (G) at (-1,.8) {$A_n$};
  \node (L) at (-1,0) {$C_n$};
  \foreach \y in {0,.8} {
  \foreach \x in {0,1,2,4,5} { \fill (\x,\y) circle (.7mm);}
  \draw (1,\y) -- (2,\y);
  \draw[dotted] (2,\y) -- (4,\y);
  \draw (4,\y) -- (5,\y);
  }
  \draw (0,.8) -- (1,.8);
  \draw (0,0.05) -- (1,0.05);
  \draw (0,-0.05) -- (1,-0.05);
  \node (0) at (0,-.6) {$i_1$}; \node (1) at (1,-.6) {$i_2$}; \node (3) at (2,-.6) {$i_3$}; \node (4) at (4,-.6) {$i_{n-1}$}; \node (5) at (5,-.6) {$i_n$}; \node (dot) at (3,-.6) {$\cdots$};
\end{tikzpicture}
\end{center}

\noindent Associated to every Coxeter matrix, there is a finitely presented group, the Coxeter group.

\begin{Definition}
    Let $I$ be a finite set and let $M$ be a Coxeter matrix. The associated \emph{Coxeter group $W$} is a finitely presented group with generator set $S=\{s_i:i\in I\}$ and with presentation
    \[
    W = \langle S \,|\, (s_is_j)^{m_{i,j}}=1\quad\forall i,j\in I \rangle.
    \]
    The pair $(W,S)$ is called a \emph{Coxeter system}, the cardinality $|I|$ is called its \emph{rank}. The set of all cosets
	\[
		\Sigma(W,S)=\{w\langle S'\rangle: w\in W, S'\subseteq S\},
	\]
    partially ordered by reverse inclusion forms a simplicial complex, the \emph{Coxeter complex} associated to $(W,S)$.

    If $\sigma\in\Sigma(W,S)$ is of the form $\sigma=w\langle \{s_i : i \in I' \}\rangle$, we write $\type(\sigma) = I\setminus I'$. In particular, the types of vertices are single elements of $I$.

    Top-dimensional simplices of $\Sigma(W,S)$ are called \emph{chambers}, codimension-1-simplices are called \emph{panels}.
\end{Definition}

\noindent The definition of opposition in spherical buildings is central to our discussion later on.

\begin{Definition}
    We say that $W$ is \emph{spherical} if $W$ is finite. In this case $\Sigma(W,S)$ can be realised naturally as a triangulated sphere of dimension $|I|-1$. Then the antipodal map of the sphere induces a simplicial involutory automorphism $\rho$ of $\Sigma(W,S)$. A simplex $\sigma$ and its image $\rho(\sigma)$ are said to be \emph{opposite}. The notion of opposition can also be expressed in a combinatorial fashion.
\end{Definition}

\noindent Buildings are simplicial complexes covered by the union of their apartments, which are copies of a fixed Coxeter complex.

\begin{Definition}
	A simplicial complex $\Delta$ together with a collection of subcomplexes $\cA$ called \emph{apartments} is a \emph{building} if
	\begin{itemize}
		\item Every apartment is a Coxeter complex.
		\item Every two simplices of $\Delta$ are contained in a common apartment.
		\item For any two apartments $A_1,A_2\in \cA$ containing a common chamber, there is an isomorphism $A_1\rightarrow A_2$ fixing $A_1\cap A_2$ pointwise.
	\end{itemize}
    The building $\Delta$ is called \emph{thick} if every panel is contained in at least three chambers. We will call a building \emph{weak} if it is not necessarily thick.

\end{Definition}

\noindent The axioms force all apartments to be isomorphic, there is in particular a unique Coxeter system $(W,S)$ associated to $\Delta$.

\begin{Definition}
    If this Coxeter group $W$ is spherical, then $\Delta$ is called spherical as well. Likewise, we say that the type of $\Delta$ is the type of $(W,S)$. The \emph{rank} of the building is the cardinality of the type set $I$.
\end{Definition}

\noindent We will need the following properties of buildings:

\begin{Lemma}
    The link of a simplex $\lk_\Delta(\sigma)$ is again a building of type $I\setminus \type(\sigma)$.
\end{Lemma}

\begin{Theorem}[Solomon-Tits]
    A spherical building of rank $n$ has the homotopy type of a bouquet of $(n-1)$-spheres.
\end{Theorem}

\noindent The following example will be discussed again in section \ref{sec:general_linear_groups}.

\begin{Example}
    The simplest example of a rank $n$ spherical building is the flag complex $\Delta$ over the $n$-dimensional projective space over any division ring $D$. Simplices in $\Delta$ are then ascending flags of subspaces of $D^{n+1}$:
    \[
    \Delta = \{ (V_1 \subsetneq V_2 \subsetneq \cdots\subsetneq V_k) : V_1\neq 0, V_k\subsetneq D^{n+1} \}
    \]
    Apartments in $\Delta$ consist of all flags whose elements can be expressed as spans of a fixed basis of $D^{n+1}$. The associated Coxeter group $W$ is the symmetric group on $n+1$ letters permuting the basis vectors. This building is of type $A_n$.

    Two simplices $(V_1\subsetneq V_2 \subsetneq\cdots\subsetneq V_k)$ and $(V'_1\subsetneq V'_2 \subsetneq\cdots\subsetneq V'_k)$ are opposite if and only if $V_i\bigoplus V'_{k+1-i} = D^{n+1}$ for all $i=1,\ldots,k$.

    The type of a simplex $\sigma=(V_1\subsetneq V_2 \subsetneq\cdots\subsetneq V_k)$ is the set of dimensions of subspaces in the flag: $\type(\sigma)=\{ i_l : \dim(V_i)=l, 1\leq i\leq k\}$.
\end{Example}

\noindent Our main result will involve groups with weak spherical Tits systems or, equivalently, groups acting strongly transitively on weak spherical buildings.

\begin{Definition}
    A group $G$ acts \emph{strongly transitively} on a weak spherical building $\Delta$ if it acts transitively on pairs consisting of an apartment and a chamber in it.
\end{Definition}
    
\noindent Strongly transitive actions give rise to a group theoretic datum called a Tits system.

\begin{Definition}
    Let $G$ be a group. Let $B$ and $N$ be subgroups of $G$ such that their intersection $H=B\cap N$ is normal in $N$. Assume also that the quotient group $W\coloneq N/H$ is generated by a set $S\subset W$. The quadruple $(G,B,N,S)$ is called a \emph{weak Tits system} if
    \begin{itemize}
        \item $G=\langle B \cup N\rangle$,
        \item $(W,S)$ is a Coxeter system and
        \item for $s\in S$ and $w\in W$, we have $BsBwB\subseteq BwB \sqcup BswB$.
    \end{itemize}
    If in addition $sBs\neq B$ for all $s\in S$, we call $(G,B,N,S)$ a \emph{Tits system} for $G$.
\end{Definition}

\noindent Given a strongly transitive action of a group $G$ on a weak spherical building, we can construct a Tits system as follows.

\begin{Construction}
    Fix an apartment $\Sigma$ and a chamber $c_0\in\Sigma$. Denote by $B$ the stabiliser of $c_0$ and by $N$ the normaliser of $\Sigma$. Their intersection $H=B\cap N$ is the pointwise stabiliser of $\Sigma$ and the group $H$ is normal in $N$. The quotient $N/H$ acts chamber-regularly on the apartment $\Sigma$ and is hence isomorphic to the associated Coxeter group $W$, where we have fixed the generating set $S$.
\end{Construction}

\begin{Theorem}
    The quadruple $(G,B,N,S)$ is a weak Tits system for $G$. If the building is thick, then $(G,B,N,S)$ is a Tits system. Conversely, given a (weak) Tits system of spherical type for $G$, a (weak) spherical building can be constructed on which $G$ acts strongly transitively. 
\end{Theorem}

\noindent The proof of this result can be found in \cite[Chapter 6]{AB:B:08}.

\begin{Remark}
    We will later describe the buildings and the weak Tits systems we use explicitly. The class of groups with weak spherical Tits systems includes general and special linear groups over division rings and unitary groups over hyperbolic modules, in particular symplectic groups and special orthogonal groups of maximal Witt index.

    Note that the latter groups have a (non-weak) Tits system of type $D_n$, but also a weak Tits system of type $C_n$, as described in \cite[6.7]{AB:B:08}.
\end{Remark}

\begin{Definition}
    Any stabiliser in $G$ of a simplex in $\Delta$ is called a \emph{parabolic subgroup}.
\end{Definition}

\begin{Example}
    The general linear group $\Gl_{n+1}(D)$ acts strongly transitively on the building $\Delta$ of type $A_n$, which is the flag complex over projective space, as discussed in the previous section. If $e_1,e_2,\ldots,e_{n+1}$ is the standard basis of $D^{n+1}$, the parabolic subgroup associated to the vertex $v=(\langle e_1,\ldots, e_k\rangle)$ is given by
    \[
    G_v = \begin{pmatrix}
        \Gl_{k}(D) & * \\ 0 & \Gl_{n+1-k}(D)
    \end{pmatrix}.
    \]
\end{Example}

\subsection{The opposition complex}

Let $n\geq 2$ and let $\Delta$ be a weak spherical building of rank $n$. Enumerate the type set $I=\{i_1,\ldots,i_n\}$ of $\Delta$ arbitrarily. In addition, we fix a group $G$ that acts strongly transitively on $\Delta$.

The basic geometry on which we will study the action of $G$ is not the building $\Delta$, but its associated opposition complex.

\begin{Definition}
	The \emph{opposition complex $O(\Delta)$} is the simplicial complex consisting of pairs of opposite simplices
	\[
		O(\Delta) \coloneq \{ \sigma= (\sigma^+,\sigma^-) \in \Delta\times \Delta : \sigma^+ \text{ opposite } \sigma^- \}
	\]
	with the induced inclusion relations. Set $\type((\sigma^+,\sigma^-)) \coloneq \type(\sigma^+)$.
\end{Definition}

\noindent The opposition complex $O(\Delta)$ is a $G$-simplicial complex since the $G$-action preserves opposition. The $G$-action is transitive on vertices of the same type of $O(\Delta)$, since it is transitive on pairs of opposite vertices of a fixed type in the building.

\begin{Example}
    The vertices of the opposition complex associated to the general linear group $\Gl_{n+1}(D)$ acting on the associated projective space over $D^{n+1}$ as in Section \ref{subsec:sph_buildings} are pairs of complementary subspaces of $D^{n+1}$.
\end{Example}

\noindent The significance of the opposition complex for this paper lies in the following theorem.

\begin{Theorem}[von Heydebreck, \cite{vH:HPC:03}, Theorem 3.1]\label{th:opposition_complex_is_spherical}
    The opposition complex of a weak spherical building $\Delta$ of rank $n$ is homotopy equivalent to a bouquet of $(n-1)$-spheres, we will also say that it is \emph{$(n-1)$-spherical}.
\end{Theorem}

\noindent We fix a set of representative vertices for the $G$-action.

\begin{Definition}\label{def:situation}
	We fix a standard apartment $\Sigma$ in $\Delta$ and a chamber $c_0\in\Sigma$. In the following, we write $v_p^+$ for the vertex of $c_0$ of type $\{i_p\}$ and we denote the corresponding opposite vertex in $\Sigma$ by $v_p^-$. We write $v_p=(v_p^+,v_p^-)\in O(\Delta)$ for the vertex in $O(\Delta)$.
\end{Definition}

\noindent For the inductive arguments to come, we investigate the structure of stabilisers.

\begin{Definition}
	Denote the stabilisers as follows:
\[
	L_p \coloneq G_{v_p} = G_{v_p^+} \cap G_{v_p^-}
\]
	These are intersections of two opposite parabolic subgroups and they are called \emph{Levi subgroups of $G$}.
\end{Definition}

\begin{Example}
    For the general linear group $\Gl_{n+1}(D)$ acting on the associated building $\Delta$ of type $A_n$, as above, the stabiliser of the vertex $v=(\langle e_1,\ldots,e_k\rangle,\langle e_{k+1},\ldots,e_{n+1}\rangle)$ in $O(\Delta)$ is the subgroup
    \[
        L_v = \begin{pmatrix}
        \Gl_{k}(D) & 0 \\ 0 & \Gl_{n+1-k}(D)
    \end{pmatrix}.
    \]
    Note that the Levi subgroup splits as a direct product of smaller general linear groups.
\end{Example}

\noindent The opposition complex commutes with the formation of links as follows.

\begin{Proposition}[von Heydebreck, \cite{vH:HPC:03}, Proposition 2.1]\label{prop:links_opposition_complex}
	For a simplex \\ $(\sigma^+,\sigma^-)\in O(\Delta)$ we have
	\[
		\lk_{O(\Delta)}((\sigma^+,\sigma^-)) \cong O(\lk_\Delta(\sigma^+)).
	\]
  This isomorphism is $G_{\{\sigma^+,\sigma^-\}}$-equivariant. In particular, the link $\lk_{O(\Delta)}((\sigma^+,\sigma^-))$ is $(n-1-\rank(\sigma^+))$-spherical.
\end{Proposition}

\subsection{A Filtration}

We construct an exact chain complex of $G$-modules associated to $O(\Delta)$. The construction is similar to the construction of cellular chains of a CW complex. The filtration by skeletons is replaced by a filtration by types.

\begin{Definition}
	For $1\leq p \leq n$ let $I_p=\{i_1,\ldots,i_p\}$. Write
	\[
		O(\Delta)_p \coloneq \{ \sigma\in O(\Delta) : \type (\sigma) \subseteq I_p\},
	\]
	this is a $G$-invariant filtration of $O(\Delta)$. We set $O(\Delta)_0\coloneq \emptyset$.
\end{Definition}

\noindent Observe that $O(\Delta)_p$ is of rank $p$ and hence of dimension $p-1$.

\begin{Definition}
    Let $v$ be a vertex in $O(\Delta)_p$. We define the \emph{filtrated residue}, \emph{link} and \emph{star} by:
	\begin{align*}
        R(v)_p &= R_{O(\Delta)}(v) \cap O(\Delta)_p = \{\sigma\in O(\Delta)_p: v\in\sigma\},\\
		\lk(v)_p &= \lk_{O(\Delta)}(v) \cap O(\Delta)_p = \{\sigma\in O(\Delta)_p: v\cup\sigma\in O(\Delta), v\cap\sigma=\emptyset\},\\
		\st(v)_p &= \lk(v)_p \sqcup R(v)_p.
	\end{align*}
\end{Definition}

\begin{Remark}
	From the definition it is obvious that
	\[
	\st(v)_p \cap O(\Delta)_{p-1} = \lk(v)_p = \lk(v)_{p-1}
	\]
	if $\type(v)=\{i_p\}$.
\end{Remark}

\begin{Proposition}\label{prop:filtered_homology}
	For $2\leq p\leq n$ we have
	\[
	H_i(O(\Delta)_p,O(\Delta)_{p-1}) \cong \bigoplus_{\type(v)=i_p} \tilde H_{i-1}(\lk(v)_{p-1}).
	\]
\end{Proposition}

\begin{Proof}
	We have
	\begin{align*}
		O(\Delta)_p \setminus O(\Delta)_{p-1} &= \{\sigma\in O(\Delta) : i_p \in \type(\sigma) \subseteq I_p \} \\
		&= \coprod_{\type(v)=i_p} \{\sigma\in O(\Delta)_p : v\in \sigma \} \\
		&= \coprod_{\type(v)=i_p} R(v)_p.
	\end{align*}
	Since $\st(v)_p\cap O(\Delta)_{p-1} = \lk(v)_{p-1}$, we obtain the following pushout diagram
  \[
    \xymatrix{
    \coprod_{\type(v)=i_p} \lk(v)_{p-1}\ar[r]\ar[d] & O(\Delta)_{p-1}\ar[d] \\
    \coprod_{\type(v)=i_p} \st(v)_p\ar[r] & O(\Delta)_p.
    }
  \]
	By excision, we obtain
	\begin{align*}
		H_i(O(\Delta)_p,O(\Delta)_{p-1}) &\cong \bigoplus_{\type(v)=i_p} H_i(\st(v)_p,\lk(v)_{p-1}) \\
		&\cong \bigoplus_{\type(v)=i_p} \tilde H_i(\st(v)_p / \lk(v)_{p-1}) \\
		&\cong \bigoplus_{\type(v)=i_p} \tilde H_{i-1}(\lk(v)_{p-1}).
	\end{align*}
    The last line follows from the fact that $\st(v)_p$ is the simplicial cone over $\lk(v)_{p-1}$.
\end{Proof}

\noindent This description of relative homology allows us to show that each filtration subcomplex of the opposition complex is also spherical.

\begin{Proposition}\label{prop:filtrated_opposition_complexes_are_spherical}
	For any $1\leq p \leq n$, the homology groups $\tilde H_i(O(\Delta)_p)$ are trivial except for $i=p-1$. The group $\tilde H_{p-1}(O(\Delta)_p)$ is $\Z$-free.
\end{Proposition}

\begin{Proof}
    Since the complexes $O(\Delta)_p$ are $(p-1)$-dimensional, their top-dimensional homology group $H_{p-1}(O(\Delta)_p)$ is automatically $\Z$-free.

	It remains to show that all other reduced homology groups vanish. We prove this by induction on $n=\rank(\Delta)$. In any case, the statement is true for $O(\Delta)_n=O(\Delta)$ by Theorem \ref{th:opposition_complex_is_spherical} and it is trivial for $O(\Delta)_1$.

	Combining these facts, we obtain the statement for $n=2$. Now assume $n\geq 3$. We prove the statement for $O(\Delta)_p$ for all $2\leq p\leq n$ by reverse induction. The case $p=n$ is already known by Theorem \ref{th:opposition_complex_is_spherical}.

	Hence assume we have the statement for $O(\Delta)_p$ and prove it for $O(\Delta)_{p-1}$. First of all note that by Proposition \ref{prop:filtered_homology}, we have
	\[
	H_i(O(\Delta)_p,O(\Delta)_{p-1}) \cong \bigoplus_{\type(v)=i_p} \tilde H_{i-1}(\lk(v)_{p-1})
	\]
	and we have the induction hypothesis for $\lk(v)_{p-1}$, since $\rank(\lk_\Delta(v^+))=\rank(\Delta)-1$.
	We see in particular that $H_i(O(\Delta)_p,O(\Delta)_{p-1})$ vanishes for $i\neq p-1$.

    By the long exact sequence for the pair $(O(\Delta)_p,O(\Delta)_{p-1})$, we obtain that $\tilde H_i(O(\Delta)_{p-1})$ vanishes for $i\neq p-2$.
\end{Proof}

\noindent The following modules $M_p$ will be the coefficient modules in the spectral sequence.

\begin{Definition}
	We define a sequence of $L_p$-modules as follows:
	\[
		M_p \coloneq\begin{cases}
			\Z & p=1 \\
			\tilde H_{p-2}(\lk(v_p)_{p-1}) & 2\leq p \leq n.
		\end{cases}
	\]
    These are $L_p$-modules, since $L_p$ stabilises $\lk(v_p)$ and is type-preserving. It hence also stabilises the subcomplex $\lk(v_p)_{p-1}$. Note additionally that $M_p$ is $\Z$-free by the previous proposition.
\end{Definition}

\noindent With this definition, we obtain a simple description of the relative homology modules.

\begin{Proposition}\label{prop:structure_of_filtrated_homology}
	For $1\leq p\leq n$, we have
	\[
	H_i(O(\Delta)_p,O(\Delta)_{p-1}) \cong\begin{cases}
		0 & i \neq p-1\\ \Z G \otimes_{L_p} M_p & i=p-1.
	\end{cases}
	\]
\end{Proposition}

\begin{Proof}
	If $p=1$, then $O(\Delta)_1 = \coprod_{\type(v)=i_1} \{v\}$. Then obviously
	\[
	H_0(O(\Delta)_1) \cong \Z G \otimes_{L_1} \Z
	\]
	since $G$ acts transitively on pairs of opposite vertices of the same type.

	For $p>1$, by Proposition \ref{prop:filtered_homology}, we have
	\begin{align*}
        H_i(O(\Delta)_p,O(\Delta)_{p-1}) &\cong \bigoplus_{\type(v)=i_p} \tilde H_{i-1}(\lk(v)_{p-1}), \\
        &\cong \Z G \otimes_{L_p} \tilde H_{i-1}(\lk(v_p)_{p-1} ),
	\end{align*}
    again since $G$ acts transitively on pairs of opposite vertices of the same type. The claim now follows from Proposition \ref{prop:filtrated_opposition_complexes_are_spherical}.
\end{Proof}

\noindent The filtration allows us to obtain an exact complex of $G$-modules, which will be used to construct a relative spectral sequence.

\begin{Definition}\label{def:exact_chain_complexes}
	Consider the following sequence of $G$-modules:
	\[
		C_p\coloneq\begin{cases}
			H_{n-1}(O(\Delta)) & p=n+1 \\
			H_{p-1}(O(\Delta)_p,O(\Delta)_{p-1}) & 1 \leq p \leq n \\
			\Z & p=0 \\
			0 & \text{otherwise.}
		\end{cases}
	\]
	We then have
	\[
	C_p \cong \Z G \otimes_{L_p}M_p
	\]
  for $1\leq p \leq n$ by Proposition \ref{prop:structure_of_filtrated_homology}.
\end{Definition}

\noindent As for cellular chains, there is a chain complex structure on $C_*$. Note that, in contrast to the situation of cellular chains, we have added the modules $C_0$ and $C_{n+1}$ to make the chain complex $C_*$ exact.

\begin{Lemma}\label{l:kdelta_exact_chain_complex}
	There is a boundary map $\partial^C_*$ on $C_*$, which makes $C_*$ into an exact chain complex of $G$-modules.
\end{Lemma}

\begin{Proof}
    The filtration $(O(\Delta)_p)_{p\in\Z}$ induces a $G$-equivariant filtration on cellular chains of $O(\Delta)$. There is hence a spectral sequence of $G$-modules
	\[
	E^1_{p,q} = H_{p+q}(O(\Delta)_{p+1},O(\Delta)_p) \quad\Rightarrow\quad H_{p+q}(O(\Delta)).
	\]
	By Proposition \ref{prop:filtered_homology} we obtain
	\[
	E^1_{p,q} = \begin{cases}
		H_p(O(\Delta)_{p+1},O(\Delta)_p) & q=0 \\
		0 & q \neq 0.
	\end{cases}
	\]
	The spectral sequence hence collapses on the second page, the differential maps on the first page and the edge homomorphisms form the long exact sequence.
\end{Proof}

\begin{Remark}
    A closer comparison to the situation of cellular chains shows that this boundary map is given by the composition
    \[
    H_p(O(\Delta)_{p+1},O(\Delta)_p) \stackrel{\delta}{\rightarrow} H_{p-1}(O(\Delta)_p) \stackrel{H(\pi)}{\rightarrow} H_{p-1}(O(\Delta)_p,O(\Delta)_{p-1})
    \]
    where $\delta$ is the connecting homomorphism of the long exact sequence associated to the pair $(O(\Delta)_{p+1},O(\Delta)_p)$ and $H(\pi)$ is induced by the projection as in the long exact sequence associated to the pair $(O(\Delta)_p,O(\Delta)_{p-1})$.
\end{Remark}

\subsection{The relative spectral sequence}

For $n\geq 2$, let $G$ be a group with a weak Tits system of rank $n+1$ with associated building $\Delta$ whose type set $I=\{i_1,\ldots,i_{n+1}\}$ is ordered arbitrarily. We adapt the notation from the previous section.

In this general setting, of course, a part of the problem of homological stability is its precise formulation. Which subgroups of a given group $G$ should be considered as the ones to yield stability? These subgroups cannot be expressed explicitly in this generality, they depend on the chosen series of groups. Consequently, we allow for a certain amount of flexibility in the choice of the subgroup $G'$.

\begin{Definition}
    Let $G'$ be any subgroup of $L_{n+1}$ that still acts strongly transitively on the link $\Delta'\coloneq\lk_\Delta(v_{n+1}^+)$, which is a building of type $I_n=\{i_1,\ldots,i_n\}$. In particular, $G'$ admits a weak Tits system of type $I_n$.
	We call the pair $(G,G')$ a \emph{stability pair}.
\end{Definition}

\begin{Remark}
  We will always choose $G'$ such that this assumption of strong transitivity is obviously fulfilled. In fact, if the group $G$ admits a root datum, then every Levi subgroup does as well by \cite[6.2.3]{Rem:GKM:02}. The group $L_{n+1}^\dagger=[L_{n+1},L_{n+1}]$ is then the associated little projective group which also admits a root datum and we will always choose $G'$ to contain $L_{n+1}^\dagger$.
\end{Remark}

\noindent The aim of this section is to associate a spectral sequence to every stability pair. It will be a relative version of the spectral sequences in \cite[VII.5]{Bro:CoG:82}.

\begin{Definition}
	By Proposition \ref{prop:links_opposition_complex}, we have
	\[
	O(\Delta)' \coloneq \lk_{O(\Delta)}(v_{n+1}) \cong O(\lk_\Delta(v_{n+1}^+))=O(\Delta').
	\]
	This isomorphism is $G'$-equivariant.
	Define a type filtration on $O(\Delta')$ analogously to the one defined on $O(\Delta)$. In addition, the filtration on $O(\Delta)$ induces a filtration on $O(\Delta)'$. The isomorphism is then also filtration-preserving.
\end{Definition}

\noindent Now we consider the exact chain complexes $C_*$ and $C'_*$ associated to the groups $G$ and $G'$ as in Definition \ref{def:exact_chain_complexes}. We obtain the following descriptions:
\begin{align*}
		C_p&\coloneq\begin{cases}
			H_n(O(\Delta)) & p=n+2 \\
			H_{p-1}(O(\Delta)_p,O(\Delta)_{p-1}) & 1 \leq p \leq n+1 \\
			\Z & p=0 \\
			0 & \text{otherwise.}
		\end{cases}\\
	\intertext{and}
		C_p'&\coloneq\begin{cases}
			H_{n-1}(O(\Delta)') & p=n+1 \\
			H_{p-1}(O(\Delta)'_p,O(\Delta)'_{p-1}) & 1 \leq p \leq n \\
			\Z & p=0 \\
			0 & \text{otherwise.}
		\end{cases}
\end{align*}
	We see that
	\[
  C_p' \cong \Z G' \otimes_{L_p'} M_p\qquad\text{ and }\qquad C_p\cong \Z G \otimes_{L_p} M_p
	\]
  for $1\leq p\leq n$ and $C'_{n+1} = H_{n-1}(\lk_{O(\Delta)}(v_{n+1})) = M_{n+1}$, where
  \[
  M_p = \tilde H_{p-2}(\lk_{O(\Delta)}(v_p)_{p-1}) \cong \tilde H_{p-2}(\lk_{O(\Delta)'}(v_p)_{p-1})\qquad\text{for } 2\leq p\leq n
  \]
  are the same modules for both groups $G$ and $G'$.

\begin{Lemma}
	The inclusion $O(\Delta)'\hookrightarrow O(\Delta)$ induces a $G'$-equivariant chain map
  \[
  \iota: C_*' \rightarrow C_*.
  \]
\end{Lemma}

\begin{Proof}
	The inclusion $O(\Delta)'\hookrightarrow O(\Delta)$ is $G'$-equivariant and filtration-preserving. The inclusions of pairs then induce maps
	\[
	\iota_p: \underbrace{H_{p-1}(O(\Delta)'_p,O(\Delta)'_{p-1})}_{C_p'} \rightarrow \underbrace{H_{p-1}(O(\Delta)_p,O(\Delta)_{p-1})}_{C_p}
	\]
	for $1\leq p \leq n$, which are compatible with the boundary maps $\partial^C$ and $\partial^{C'}$. Of course $\iota_p=0$ for $p\leq -1$ and $p\geq n+2$. Consider
	\[
	\xymatrix{
		0 & \Z\ar[l] & H_0(O(\Delta)_1) \ar[l]_-{\partial^C_1} & \ar[l]\cdots \\
		0 & \Z\ar[l]\ar@.[u]^{\iota_0} & H_0(O(\Delta)'_1)\ar[l]_-{\partial^{C'}_1}\ar[u]^{\iota_1} & \ar[l]\cdots
	}
	\]
  Here $\iota_0\coloneq \partial^C_1 \circ \iota_1 \circ (\partial^{C'}_1)^{-1}$ is well-defined by a diagram chase. For $p=n+1$, consider
	\[
		\xymatrix{
        \cdots & H_{n-1}(O(\Delta)_n,O(\Delta)_{n-1})\ar[l] & H_n(O(\Delta)_{n+1},O(\Delta)_n)\ar[l]_-{\partial^C_{n+1}} & H_{n}(O(\Delta)_{n+1})\ar[l] \\
        \cdots & H_{n-1}(O(\Delta)'_n,O(\Delta)'_{n-1})\ar[l]\ar[u]^{\iota_n} & \tilde H_{n-1}(O(\Delta)'_n) \ar[l]_-{\partial^{C'}_{n+1}}\ar@.[u]_{\iota_{n+1}} & 0.\ar[l]
		}
	\]
    Note that $O(\Delta)'_n=\lk_{O(\Delta)}(v_{n+1}) = \lk(v_{n+1})_n$, so
    \[
    \tilde H_{n-1}(O(\Delta)'_n) \cong H_n\bigl(\st(v_{n+1})_{n+1},\lk(v_{n+1})_n\bigr)
    \]
    as in the proof of Proposition \ref{prop:filtered_homology}. The inclusion of pairs 
    \[
        \bigl(\st(v_{n+1})_{n+1},\lk(v_{n+1})_n\bigr) \hookrightarrow \bigl(O(\Delta)_{n+1},O(\Delta)_n\bigr)
    \]
    then induces the map $\iota_{n+1}$ on homology. A closer inspection using the explicit description of the boundary maps above shows that $\iota$ is a chain map.
  \end{Proof}

\begin{Remark}
    It is not difficult to see that, under the above isomorphisms, the map $\iota_p$ is induced by the canonical inclusions 
    \[
    \Z G'\otimes_{L'_p} M_p \hookrightarrow \Z G \otimes_{L_p} M_p
    \]
    for $1\leq p \leq n$ and that $\iota_{n+1}: M_{n+1} \rightarrow \Z G \otimes_{L_{n+1}} M_{n+1}$ is given by $m\mapsto 1\otimes_{L_{n+1}} m$.
\end{Remark}

\noindent With these chain complexes and the chain map, we are able to construct the following spectral sequence.

\begin{Theorem}[Stability pair spectral sequence]\label{th:stability_pair_spectral_sequence}
    For each stability pair $(G,G')$, there is a first-quadrant spectral sequence
	\[
	E^1_{p,q}=\begin{cases}
		H_q(G,G';\Z) & p=0 \\
		H_q(L_p,L_p';M_p) & 1\leq p\leq n \\
    H_q(L_{n+1},G';M_{n+1}) & p=n+1
	\end{cases}
    \]
    which converges to zero.
\end{Theorem}

\begin{Proof}
    Let $F_*(G)$ and $F_*(G')$ be the standard resolutions of $\Z$ over $\Z G$ and $\Z G'$, respectively. Consider the two double complexes $F_*(G')\otimes_{G'} C'_*$ and $F_*(G)\otimes_G C_*$ and the map of double complexes
    \begin{align*}
    i: F(G')\otimes_{G'} C' &\rightarrow F(G)\otimes_G C \\
    f \otimes_{G'} c &\mapsto f \otimes_G \iota(c)
\end{align*}
    induced by the inclusion $G'\hookrightarrow G$ and by $\iota$.
    We can apply Theorem \ref{th:relative_spectral_sequence} to obtain the relative spectral sequence
    \[
    E^1_{p,q} = H_q(\Cone_*(i_{\bullet,p})) \Rightarrow 0,
    \]
    converging to zero since $C$ and $C'$ are exact chain complexes. All that remains is the description of the first page terms. For $1\leq p\leq n$, we apply \cite[1.5.8]{Wei:IHA:94} to obtain
    \[
    E^1_{p,q} \cong H_q\bigl( (F_*(G) \otimes_G C_p ) / (F_*(G') \otimes_{G'} C'_p)\bigr) \cong H_q\bigl( (F_*(G) \otimes_{L_p} M_p ) / (F_*(G') \otimes_{L'_p} M_p ) \bigr),
    \]
    which is isomorphic to $H_q(L_p,L'_p;M_p)$ as in the proof of Theorem \ref{th:relative_lyndon_hochschild_serre}. For $p=n+1$, note that
    \begin{align*}
    E^1_{n+1,q} &= H_q\bigl( (F_*(G) \otimes_G C_{n+1} ) / (F_*(G') \otimes_{G'} C'_{n+1})\bigr) \\
    &\cong H_q\bigl( (F_*(G) \otimes_{L_{n+1}} M_{n+1} ) / (F_*(G') \otimes_{G'} M_{n+1} ) \bigr),
\end{align*}
    which is isomorphic to $H_q(L_{n+1},G';M_{n+1})$.
\end{Proof}

\begin{Remark}\label{rem:e1_vanishes_at_n+1}
  Since relative $H_0$ vanishes always, we have $E^1_{p,0}=0$ for all $0\leq p\leq n+1$.
\end{Remark}

\section{Group theory}

In the third part, we will apply the results of the first two parts to specific series of groups. We will only sketch the original application to general linear groups by Charney in \cite{Cha:HSD:80}. Instead, we cite a strong stability theorem for general linear groups by Sah, which we will use for the following proofs. In the following sections, we then prove homological stability for special linear groups, for unitary groups and for special orthogonal groups.

\subsection{General linear groups}\label{sec:general_linear_groups}

Let $D$ be any division ring. We consider the case where $G=\Gl_{n+2}(D)$ for $n\geq 2$. The associated building $\Delta$ is the flag complex over $(n+1)$-dimensional projective space, the opposition complex $O(\Delta)$ consists of pairs of complementary vector subspaces of $D^{n+2}$. We choose a basis $e_1,e_2,\ldots,e_{n+2}$ of $D^{n+2}$ corresponding to the standard apartment $\Sigma$. We choose the type filtration and the chamber $c_0$ in $\Sigma$ with vertices $v_p^+ = \langle e_1,\ldots,e_p \rangle$, hence $v_p^- = \langle e_{p+1},\ldots,e_{n+2}\rangle$. Consequently
\[
	L_p = \begin{pmatrix}
		\Gl_p(D) & 0 \\ 0 & \Gl_{n+2-p}(D)
	\end{pmatrix}
\] for $1\leq p \leq n+1$. In particular \[ L_{n+1}\cong \begin{pmatrix} \Gl_{n+1}(D) & 0 \\ 0 & D^\times
\end{pmatrix}.
\]
We choose $G' = \begin{pmatrix}
	\Gl_{n+1}(D) & 0 \\ 0 & 1
\end{pmatrix}$ which implies that
\[
	L'_p = \begin{pmatrix}
		\Gl_p(D) & 0 & 0 \\
		0 & \Gl_{n+1-p}(D) & 0 \\
	  0 & 0 & 1 \\
	\end{pmatrix}.
\]
Obviously, we obtain $L_p\cong \Gl_p(D) \times \Gl_{n+2-p}(D)$ and $L'_p \cong \Gl_p(D)\times \Gl_{n+1-p}(D)$. By Theorem \ref{th:stability_pair_spectral_sequence} we obtain as in \cite{Cha:HSD:80} the following spectral sequence.

\begin{Theorem}[Charney, Theorem 2.2 in \cite{Cha:HSD:80}]
	For any $n\geq 2$, there is a first-quadrant spectral sequence
	\[
	E^1_{p,q} =\begin{cases}
		H_q\bigl(\Gl_{n+2}(D),\Gl_{n+1}(D);\Z\bigr) & p=0 \\
		H_q\bigl(\Gl_{p}(D)\times \Gl_{n+2-p}(D),\Gl_p(D)\times \Gl_{n+1-p}(D); M_p\bigr) & 1\leq p \leq n\\
        H_q\bigl(\Gl_{n+1}(D)\times D^\times, \Gl_{n+1}(D);M_{n+1}\bigr) & p=n+1
	\end{cases}
	\]
	which converges to zero.
\end{Theorem}

\begin{Remark}
    It is not difficult to see that the maps
	\[
  E^1_{1,q} = H_q\bigl(D^\times \times \Gl_{n+1}(D), D^\times \times \Gl_n(D);\Z\bigr) \rightarrow H_q\bigl(\Gl_{n+2}(D),\Gl_{n+1}(D);\Z\bigr) = E^1_{0,q}
	\]
	are induced by the inclusion of pairs, which is an important ingredient in the homological stability proof by Charney described below.
\end{Remark}

\noindent Note that the vertices of $\lk(v_p)_{p-1}$ are given by
\[
\lk(v_p)_{p-1}^0 = \bigl\{ (V,W) : V\subsetneq \langle e_1,\ldots,e_p\rangle, W \supsetneq \langle e_{p+1},\ldots,e_{n+2}\rangle, V\oplus W=D^{n+2}\bigr\}.
\]
In particular, the following factors of the Levi subgroups
\[ \begin{pmatrix}
	\one_p & 0 \\ 0 & \Gl_{n+2-p}(D)
\end{pmatrix}\qquad\text{and}\qquad\begin{pmatrix}
	\one_p & 0 & 0\\ 0 & \Gl_{n+1-p}(D) & 0 \\ 0 & 0 & 1
\end{pmatrix}
\]
act trivially on $\lk(v_p)_{p-1}$ and hence on $M_p$. We can hence apply the relative Lyndon\slash Hochschild-Serre spectral sequence (Theorem \ref{th:relative_lyndon_hochschild_serre}) to obtain a description of the first page of the aforementioned spectral sequence in terms of integral relative group homology of smaller general linear groups. This can be used to apply an ingenious ``bootstrap procedure'' to prove homological stability inductively for general linear groups as in \cite{Cha:HSD:80}. The result, originally proved for Dedekind rings, is

\begin{Theorem}[Charney, Theorem 3.2 in \cite{Cha:HSD:80}]
	For $n\geq 3k$ we have
	\[
	H_k(\Gl_{n+1}(D),\Gl_{n}(D);\Z) = 0.
	\]
\end{Theorem}

\noindent This theorem is far from optimal, to the author's knowledge, the following result by Sah is the strongest for division rings with infinite centre.

\begin{Theorem}[Sah, Appendix B in \cite{Sah:HcL:86}]\label{th:sah}
	If $D$ has infinite centre then $n\geq k$ implies
	\[
	H_k(\Gl_{n+1}(D),\Gl_n(D);\Z) =0.
	\]
\end{Theorem}
The discussion in \cite[Appendix B]{Sah:HcL:86} is very brief. A detailed exposition of this result can be found in \cite{Ess:HSG:06}. A different proof of this result can be found in \cite[2.3]{Knu:HLG:01}. For arbitrary division rings, compare with this very special case of a general result by van der Kallen:

\begin{Theorem}[van der Kallen, \cite{vdK:HSL:80}]\label{th:vdk_gln}
  For any division ring $D$ we have
	\[
	H_k(\Gl_{n+1}(D),\Gl_n(D);\Z) =0.
	\]
  for $n\geq 2k$.
\end{Theorem}

\subsection{Special linear groups}

In this section, we consider special linear groups $\Sl_{n+2}(D)$ over infinite fields $D$. It is well known that the groups $\Sl_{n+2}(D)$ also act strongly transitively on the associated building $\Delta$ from the previous section. The opposition complex $O(\Delta)$ and hence also the modules $M_p$ are the same as for general linear groups. In the case $G=\Sl_{n+2}(D)$ for $n\geq 2$, the Levi subgroups admit the following structure:
\[
L_p = \biggl\{\begin{pmatrix}
	A & 0 \\ 0 & B
\end{pmatrix} : A \in \Gl_{p}(D), B \in \Gl_{n+2-p}(D), \det(A)\det(B)=1 \biggr\},
\]
which splits as a semidirect product
\begin{align*}
L_p &= \biggl\{\begin{pmatrix}
  A & 0 & 0 \\ 0 & \det(A)^{-1}  & 0 \\ 0 & 0 & \one_{n+1-p}
\end{pmatrix} : A \in \Gl_p(D) \biggr\}\ltimes\begin{pmatrix}
	\one_{p} & 0 \\ 0 & \Sl_{n+2-p}(D)
      \end{pmatrix}\\
		& \cong \Gl_{p}(D)  \ltimes \Sl_{n+2-p}(D),
\end{align*}
and the group $\Sl_{n+2-p}(D)$ acts trivially on $M_p$ by the same argument as in Section \ref{sec:general_linear_groups}. In particular
\[
L_{n+1} = \biggl\{\begin{pmatrix}
  A & 0 \\ 0 & \det(A)^{-1}
\end{pmatrix} : A\in\Gl_{n+1}(D) \biggr\}\text{, and we choose } G'=\begin{pmatrix}
  \Sl_{n+1}(D) & 0 \\ 0 & 1
\end{pmatrix}.
\]
We obtain
\begin{align*}
L_{p}' &=\biggl\{\begin{pmatrix}
  A & 0 & 0 \\ 0 & \det(A)^{-1} & 0 \\ 0 & 0 & \one_{n+1-p}
    \end{pmatrix} : A \in \Gl_p(D) \biggr\}\ltimes\begin{pmatrix}
	\one_{p} & 0 &0 \\ 0 & \Sl_{n+1-p}(D) &0 \\ 0 & 0 & 1
  	\end{pmatrix}\\
		&\cong \Gl_{p}(D) \ltimes \Sl_{n+1-p}(D),
	\end{align*}
again with $\Sl_{n+1-p}(D)$ acting trivially on $M_p$. In particular, we have
\begin{align*}
  L_1 &= \biggl\{\begin{pmatrix}
    \det(A)^{-1} & 0 \\ 0 & A
  \end{pmatrix} : A\in\Gl_{n+1}(D) \biggr\}\cong \Gl_{n+1}(D) \\
  L'_1&=\biggl\{\begin{pmatrix}
    \det(A)^{-1} & 0 & 0 \\ 0 & A & 0 \\ 0 & 0 & 1
  \end{pmatrix} : A \in \Gl_{n}(D) \biggr\}\cong \Gl_{n}(D).
\end{align*}

\noindent As in \cite{Cha:HSD:80}, we apply Theorem \ref{th:stability_pair_spectral_sequence} to obtain

\begin{Theorem}[Charney, Theorem 2.3 in \cite{Cha:HSD:80}]\label{th:sln_spectral_sequence}
	For any $n\geq 2$, there is a first-quadrant spectral sequence
	\[
	E^1_{p,q} =\begin{cases}
		H_q\bigl(\Sl_{n+2}(D),\Sl_{n+1}(D);\Z\bigr) & p=0 \\
		H_q\bigl(\Gl_{p}(D)\ltimes \Sl_{n+2-p}(D),\Gl_p(D)\ltimes \Sl_{n+1-p}(D); M_p\bigr) & 1\leq p \leq n\\
    H_q\bigl(\Gl_{n+1}(D),\Sl_{n+1}(D);M_{n+1}\bigr) & p=n+1
	\end{cases}
	\]
	which converges to zero.
\end{Theorem}

\noindent Charney then uses this spectral sequence and a version of the relative Lyndon\slash Hochschild-Serre spectral sequence (Theorem \ref{th:relative_lyndon_hochschild_serre}) to prove homological stability for $n\geq 3k$. One can improve this result using Theorem \ref{th:sah}, however.

\begin{Theorem}\label{th:sln_result}
  If $D$ is an infinite field, then $n\geq 2k-1$ implies
  \[
	  H_k(\Sl_{n+1}(D),\Sl_{n}(D);\Z)=0.
  \]
\end{Theorem}

\begin{Proof}
  In spite of the formulation of the theorem we will prove equivalently that $n\geq 2k-2$ implies
  \[
  H_k(\Sl_{n+2}(D),\Sl_{n+1}(D);\Z)=0.
  \]
  Since relative $H_0$ vanishes always and since $H_1(\Sl_{n+2}(D);\Z)=0$ for all $n\geq 0$ by \cite[2.2.3]{HOM:CGK:89}, we can start an induction over $k$.

  Let $k\geq 2$. As induction hypothesis, assume that
  \[
  H_l(\Sl_{n+2}(D),\Sl_{n+1}(D);\Z)=0\quad\text{ for all $l<k$ and all $n\geq 2l-2$.}
  \]
  We will show that $H_k(\Sl_{n+2}(D),\Sl_{n+1}(D);\Z)=0$ for all $n\geq 2k-2$. For any $n\geq 2k-2\geq 2$, by Theorem \ref{th:sln_spectral_sequence}, we have the spectral sequence $E^1_{p,q}$ converging to zero. Note first of all that
  \[
  E^1_{0,q} = H_q(\Sl_{n+2}(D),\Sl_{n+1}(D);\Z).
  \]
  In particular, we have to prove that $E^1_{0,k}=0$. Since the spectral sequence converges to zero, it is enough to prove $E^1_{p,q}=0$ when $p+q\leq k+1$ and $p\geq 1$, which will then imply $E^1_{0,q}=0$ for all $0\leq q\leq k$. The situation is illustrated in Figure \ref{fig:sln}. We will prove the vanishing of the modules separately for the regions I, II and III.

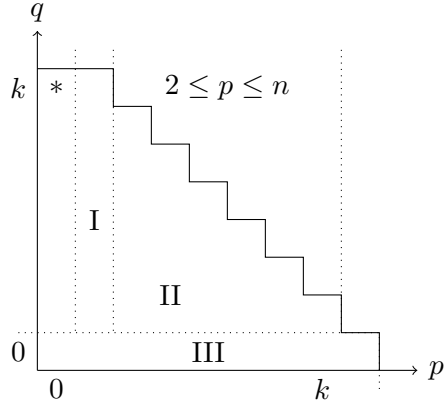
\begin{figure}
\centering
    \begin{tikzpicture}[scale=0.5,font=\small]
    \draw[->] (0,0) -- (0,9);
    \draw[->] (0,0) -- (10,0);
    \draw (0,8) -- (2,8) -- (2,7) -- (3,7) -- (3,6) -- (4,6) -- (4,5) -- (5,5) -- (5,4) -- (6,4) -- (6,3) -- (7,3) -- (7,2) -- (8,2) -- (8,1) -- (9,1) -- (9,0);
    \node (h1) at (-.5,.5) {$0$};
    \node (hk) at (-.5,7.5) {$k$};
    \node (v1) at (.5,-.5) {$0$};
    \node (vk) at (7.5,-.5) {$k$};
    \node (p) at (10.5,0) {$p$};
    \node (q) at (0,9.5) {$q$};
    \node (leq) at (5,7.5) {$2\leq p \leq n$};
    \draw[dotted] (1,8.5) -- (1,1);
    \draw[dotted] (2,8.5) -- (2,1);
    \draw[dotted] (8,8.5) -- (8,1);
    \draw[dotted] (-.5,1) -- (9,1) -- (9,-.5);
    \node (1) at (.5,7.5) {$*$};
    \node (2) at (1.5,4) {I};
    \node (3) at (3.5,2) {II};
    \node (4) at (4.5,.5) {III};
  \end{tikzpicture}\caption{$E^1_{p,q}$ in the proof of Theorem \ref{th:sln_result}}\label{fig:sln}
\end{figure}
  For region I, that is for $p=1$ and $1\leq q\leq k$, note that
  \[
  E^1_{1,q} \cong H_q(\Gl_{n+1}(D),\Gl_{n}(D);\Z).
  \]
  Since $n\geq 2k-2\geq k\geq q$, we know that $E^1_{1,q}=0$ by Theorem \ref{th:sah}.

  For region II, that is $p+q\leq k+1$, $q\geq 1$ and $2\leq p\leq n$, we have in particular $q\leq k-1$. Additionally
  \[
    2q-2+p = (p+q) + q -2 \leq 2k-2 \leq n,
  \]
  so $n-p\geq 2q-2$. By the induction hypothesis, we hence have
  \begin{equation}
      H_q(\Sl_{n-p+2}(D), \Sl_{n-p+1}(D);\Z) = 0 \label{eq:vanishing_terms_of_relative_lhs}
  \end{equation}
  for all $p,q$ with $p+q\leq k+1$ and $p\geq 2$.
  
  For the terms $E^1_{p,q}=H_q(L_p,L'_p;M_p)$ in region II, consider the relative Lyndon\slash Hochschild-Serre spectral sequence from Theorem \ref{th:relative_lyndon_hochschild_serre}:
  \[
  L^2_{i,j} = H_i(\Gl_p(D);H_j(\Sl_{n+2-p}(D),\Sl_{n+1-p}(D);\Z)\otimes_\Z M_p) \Rightarrow H_{i+j}(L_p,L'_p;M_p).
  \]
  Note that we proved in \eqref{eq:vanishing_terms_of_relative_lhs} that $L^2_{i,j}=0$ for $j\leq q\leq k+1-p$, hence $H_q(L_p,L'_p;M_p)=0$ for all $p$, $q$ in region II. Hence $E^1_{p,q}=0$ for region II.

  Finally, note that region III, that is, all $E^1_{p,0}$ with $0\leq p\leq k+1\leq n+1$ vanish by the remark after Theorem \ref{th:stability_pair_spectral_sequence}.

  Inspecting the spectral sequence $E^1_{p,q}$ yields $E^1_{0,k}=E^{\infty}_{0,k}=0$, which proves the theorem.
\end{Proof}

\begin{Remark}
The best known stability result by Hutchinson and Tao in \cite{HT:HSS:08} for special linear groups improves this to $n\geq k$ in the case where $D$ is a field of characteristic zero. Our result above improves the known stability range for all other infinite fields by one, however. Before now, the best result known to the author is due to van der Kallen \cite{vdK:HSL:80} proving stability for $n\geq 2k$ for rings with stable rank 1.
\end{Remark}

\subsection{Unitary groups}

In this section, let $D$ be a division ring. Let $J:D\rightarrow D$ be an involution and let $\varepsilon=\pm 1$. Let $V$ be a finite-dimensional right $D$-vector space. Choose an $(\varepsilon,J)$-hermitian form $h:V\times V\rightarrow D$, that means $h$ is bi-additive and
\begin{align*}
    h(av,w) &= a^Jh(v,w),& h(v,aw)&=h(v,w)a,& h(w,v)&= h(v,w)^J\varepsilon.
\end{align*}
If the characteristic of $D$ is two, we have to assume that $h$ is trace-valued, see \cite[Chapter 6]{HOM:CGK:89} for further details. A subspace $X\leq V$ is \emph{totally isotropic} if $h(v,w)=0$ for all $v, w\in X$. Denote by $(n+1)$ the \emph{Witt index of $h$}, the maximal dimension of a totally isotropic subspace of $V$. Then $2(n+1)\leq \dim(V)$. It is known that $V$ splits non-canonically as an orthogonal sum $V=\cH_{n+1}\perp W$, where $W$ is \emph{anisotropic}, that is $h(v,v)\neq 0$ for $v\in W\setminus\{0\}$, and where $\cH_{n+1}$ is a \emph{hyperbolic module}, in particular $\dim(\cH_{n+1})=2(n+1)$. A good reference for unitary groups over division rings is \cite[Chapter 6]{HOM:CGK:89}.

\begin{Definition}
  The \emph{unitary group associated to $h$} is
  \[
  \U(V) = \{ A\in\Gl(V) : h(Av,Aw) = h(v,w) \text{ for all } v,w\in V\},
  \]
  the subgroup of $h$-preserving linear automorphisms of $V$.
\end{Definition}

\noindent It is well known that $\U(V)$ admits a Tits system of type $C_{n+1}$ which is weak if $J=\id$, $\dim(V)=2(n+1)$ and $\varepsilon\neq -1$, this is the case of orthogonal groups. We also allow this case explicitly. The corresponding building $\Delta$ is isomorphic to the flag complex over all totally isotropic subspaces of $V$. The opposition complex $O(\Delta)$ is then isomorphic to the flag complex over pairs of opposite totally isotropic subspaces of $V$, where $X\leq V$ is \emph{opposite} $Y\leq V$ if $X\oplus Y^\perp = V$.

It is known that there is a basis
\[
  e_{-(n+1)},\ldots,e_{-1},e_1,\ldots,e_{n+1}\]
  of the hyperbolic module $\cH_{n+1}$ such that
\[
  h(e_i,e_j) =\begin{cases}
    \varepsilon\qquad & \text{if } i+j=0 \text{ and } i>0 \\
    1 & \text{if } i+j=0\text{ and } i<0 \\
    0 & \text{otherwise.}
  \end{cases}
\]
This basis determines an apartment of the building $\Delta$. We choose the standard chamber $c_0$ of $\Delta$ and the type enumeration such that
\[
v_p = ( \langle e_{-(n+1)},\ldots,e_{-p}\rangle, \langle e_{p},\ldots,e_{n+1}\rangle )
\]
and we write $\cH_p=\langle e_{-p},\ldots,e_{-1},e_1,\ldots,e_p\rangle$ and $\cH_0=\{0\}$. It is a simple calculation to see that the Levi subgroups $L_p$ then admit the following structure
\begin{align*}
L_p &=\biggl\{\begin{pmatrix}
    S & 0 & 0 & 0 \\
    0 & A & 0 & B \\
    0 & 0 & S^{-J} & 0 \\
    0 & C & 0 & D
  \end{pmatrix} : S\in \Gl_{n+2-p}(D), \begin{pmatrix}
    A & B \\
    C & D
  \end{pmatrix}\in \U(\cH_{p-1}\perp W)\biggr\}\\
  &\cong \Gl_{n+2-p}(D)\times \U(\cH_{p-1}\perp W).
\end{align*}
In particular, we have
\[
L_{n+1} = \biggl\{\begin{pmatrix}
    s & 0 & 0 & 0 \\
    0 & A & 0 & B \\
    0 & 0 & s^{-J} & 0 \\
    0 & C & 0 & D
\end{pmatrix} : s\in D^\times, \begin{pmatrix}
    A & B \\
    C & D
  \end{pmatrix}\in \U(\cH_n\perp W)\biggr\}.
\]
We choose $G'$ to be
\[
  G' = \biggl\{\begin{pmatrix}
    1 & 0 & 0 & 0 \\
    0 & A & 0 & B \\
    0 & 0 & 1 & 0 \\
    0 & C & 0 & D
\end{pmatrix} : \begin{pmatrix}
    A & B \\
    C & D
  \end{pmatrix}\in \U(\cH_n\perp W)\biggr\}\cong \U(\cH_n\perp W),
\]
hence
\begin{align*}
L'_p &= \biggl\{\begin{pmatrix}
  1 & 0 & 0 & 0 & 0 & 0 \\
  0 & S & 0 & 0 & 0 & 0 \\
  0 & 0 & A & 0 & 0 & B \\
  0 & 0 & 0 & S^{-J} & 0 & 0 \\
  0 & 0 & 0 & 0 & 1 & 0 \\
  0 & 0 & C & 0 & 0 & D
\end{pmatrix} : S \in\Gl_{n+1-p}(D), \begin{pmatrix}
    A & B \\
    C & D
  \end{pmatrix}\in \U(\cH_{p-1}\perp W)\biggr\}\\
&\cong \Gl_{n+1-p}(D) \times \U(\cH_{p-1}\perp W).
\end{align*}
For $n\geq 2$, we apply Theorem \ref{th:stability_pair_spectral_sequence} to obtain a spectral sequence
\[
E^1_{p,q} \cong\begin{cases}
  H_q\bigl(\U(\cH_{n+1}\perp W),\U(\cH_{n}\perp W);\Z\bigr) & p=0 \\
  H_q\bigl(\Gl_{n+2-p}(D)\times \U(\cH_{p-1}\perp W),\\ \hspace{4cm} \Gl_{n+1-p}(D) \times \U(\cH_{p-1}\perp W); M_p\bigr) \hspace{1cm} & 1\leq p \leq n\\
  H_q\bigl(D^\times\times \U(\cH_n\perp W), \U(\cH_n\perp W);M_{n+1}\bigr) & p=n+1
\end{cases}
\]
which converges to zero. Using this spectral sequence, we obtain a generalisation of the theorems by Vogtmann for orthogonal groups in \cite{Vog:HSO:79} and \cite{Vog:SPH:81}, which were later generalised to Dedekind rings by Charney in \cite{Cha:gtV:87}.

\begin{Theorem}\label{th:un_result}
  For a division ring $D$ with infinite centre, the relative homology modules
  \[
    H_k\bigl(\U(\cH_{n+1}\perp W), \U(\cH_n\perp W);\Z\bigr)
  \]
  vanish for $n\geq 2$ if $k=1$ and for $n\geq k\geq 2$. If the centre of $D$ is finite, relative homology vanishes for $n\geq 2k$.
\end{Theorem}

\begin{Proof}
    Note first that all assumptions imply $n\geq 2$, so we can construct the spectral sequence $E^1_{p,q}$ as above. For $2\leq p\leq n$, the vertices of $\lk(v_p)_{p-1}$ are given by
  \[
  \lk(v_p)_{p-1}^0 = \bigl\{ (X,Y) : \langle e_{-(n+1)},\ldots,e_{-p}\rangle \subsetneq X, Y \supsetneq \langle e_p,\ldots,e_{n+1}\rangle, X\oplus Y^{\perp}=V\bigr\}
  \]
  In particular, the general linear group factors of $L_p$ and $L'_p$ act trivially on this filtrated link and hence on $M_p$. For $p\geq 1$, we can hence apply the relative Lyndon\slash Hochschild-Serre spectral sequence (Theorem \ref{th:relative_lyndon_hochschild_serre}) to obtain
  \[
  L^2_{i,j} = H_i( \U(\cH_{p-1}\perp W); H_j(\Gl_{n+2-p}(D),\Gl_{n+1-p}(D);\Z)\otimes_\Z M_p) \Rightarrow H_{i+j}(L_p,L'_{p};M_p).
  \]

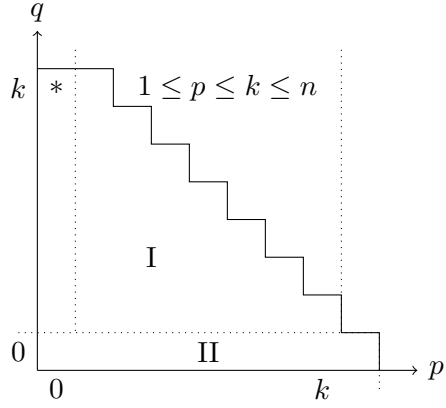
\begin{figure}
\centering
    \begin{tikzpicture}[scale=0.5,font=\small]
    \draw[->] (0,0) -- (0,9);
    \draw[->] (0,0) -- (10,0);
    \draw (0,8) -- (2,8) -- (2,7) -- (3,7) -- (3,6) -- (4,6) -- (4,5) -- (5,5) -- (5,4) -- (6,4) -- (6,3) -- (7,3) -- (7,2) -- (8,2) -- (8,1) -- (9,1) -- (9,0);
    \node (h1) at (-.5,.5) {$0$};
    \node (hk) at (-.5,7.5) {$k$};
    \node (v1) at (.5,-.5) {$0$};
    \node (vk) at (7.5,-.5) {$k$};
    \node (p) at (10.5,0) {$p$};
    \node (q) at (0,9.5) {$q$};
    \node (leq) at (5,7.5) {$1\leq p \leq k\leq n$};
    \draw[dotted] (1,8.5) -- (1,1);
    \draw[dotted] (8,8.5) -- (8,1);
    \draw[dotted] (-.5,1) -- (9,1) -- (9,-.5);
    \node (1) at (.5,7.5) {$*$};
    \node (3) at (3,3) {I};
    \node (4) at (4.5,.5) {II};
  \end{tikzpicture}\caption{$E^1_{p,q}$ in the proof of Theorem \ref{th:un_result}}\label{fig:un}
\end{figure}

Again, we want to show that $E^1_{0,k}=0$. We do this by showing that $E^1_{p,q}=0$ when $p+q\leq k+1$ and $p\geq 1$. The situation is illustrated in Figure \ref{fig:un}. We will prove the vanishing of these modules separately for the regions I and II.

First of all, consider region II, which consists of the modules $E^1_{p,0}$ with
\[
1\leq p\leq k+1\leq n+1.
\]
These modules vanish by the remark after Theorem \ref{th:stability_pair_spectral_sequence}.

Region I is given by $p+q\leq k+1$, $q\geq 1$ and $1\leq p\leq k\leq n$. We distinguish cases.

If the centre of $D$ is infinite, by Sah's result (Theorem \ref{th:sah}), we have 
\[
    H_j(\Gl_{n+2-p}(D),\Gl_{n+1-p}(D);\Z)=0
\]
for $j\leq n+1-p$. So $E^1_{p,q}$ vanishes for $p+q\leq n+1$ and $1\leq p\leq n$ by the relative Lyndon\slash Hochschild-Serre spectral sequence. Since $n\geq k$, the modules $E^1_{p,q}$ vanish in region I.

  If the centre of $(D)$ is finite, van der Kallen's result (Theorem \ref{th:vdk_gln}) implies $E^1_{p,q}=0$ for $p+2q\leq n+1$. By hypothesis, $n\geq 2k$, so $p+q\leq k+1$ and $q\leq k$ imply $p+2q\leq 2k+1\leq n+1$, which shows again that $E^1_{p,q}=0$ in region I. This proves the theorem.
\end{Proof}

\begin{Remark}
  Mirzaii and van der Kallen proved homological stability for unitary groups over local rings with infinite residue fields in \cite{MaB:HSU:02} and \cite{Mir:HSU:05} with a slightly weaker stability range of $n\geq k+1$. We restrict ourselves to division rings with infinite centre, but on the other hand we allow an anisotropic kernel.
\end{Remark}

\noindent Note that, if $W=\{0\}$ and $J=\id$, which forces $D$ to be a field, we obtain stability results for the following two special cases.

\begin{Theorem}
  For the \emph{symplectic} and \emph{orthogonal groups over an infinite field $D$}, we obtain
  \begin{align*}
    H_k\bigl(\Sp_{2n+2}(D),\Sp_{2n}(D);\Z\bigr) &=0 \\
    H_k\bigl(\Orth_{n+1,n+1}(D),\Orth_{n,n}(D);\Z\bigr) &=0
  \end{align*}
  if $k=1$ and $n\geq 2$ or $n\geq k\geq 2$. If $D$ is a finite field, then the relative homology groups vanish for $k\geq 2n$.
\end{Theorem}

\noindent By a combination of the methods for Theorem \ref{th:sln_result} and Theorem \ref{th:un_result}, results on special unitary groups can also be obtained. The following result on special orthogonal groups is particularly interesting:

\begin{Theorem}\label{th:son_result}
  For an infinite field $D$, we have
  \[
    H_k\bigl(\SO_{n+1,n+1}(D),\SO_{n,n}(D);\Z\bigr) =0
  \]
  for $k=1$ and $n\geq 2$ or $n \geq k\geq 2$. If $D$ is a finite field, then the relative homology groups vanish for $n\geq 2k$.
\end{Theorem}

\begin{Proof}
  In this case, note that $J=\id$. The Levi subgroups hence have the following simple structure:
  \[
  L_p =\biggl\{\begin{pmatrix}
    S & 0 & 0 & 0 \\
    0 & A & 0 & B \\
    0 & 0 & S^{-1} & 0 \\
    0 & C & 0 & D
  \end{pmatrix} : S\in \Gl_{n+2-p}(D), \begin{pmatrix}
    A & B \\
    C & D
  \end{pmatrix}\in \SO_{p-1,p-1}(D)\biggr\}
  \]
  Since by assumption $n\geq 2$, we can apply Theorem \ref{th:stability_pair_spectral_sequence} to obtain a relative spectral sequence converging to zero. We have to prove $E^1_{0,k}=0$.
  As in the proof of Theorem \ref{th:un_result}, the general linear factors of $L_p$, $L'_p$ act trivially on $M_p$ and we consider the relative Lyndon\slash Hochschild-Serre spectral sequence to obtain
  \[
  L^2_{i,j} = H_i\bigl( \SO_{p-1,p-1}(D); H_j(\Gl_{n+2-p}(D),\Gl_{n+1-p}(D);\Z)\otimes_\Z M_p\bigr) \Rightarrow H_{i+j}(L_p,L'_{p};M_p).
  \]
  Now, for infinite fields, note that by Theorem \ref{th:sah}, for $j\leq n+1-p$, we have $L^2_{i,j}=0$. In particular, $p+q\leq k+1\leq n+1$ and $p\geq 1$ implies $E^1_{p,q}=0$, except for $p=n+1=k+1$, where we apply the remark after Theorem \ref{th:stability_pair_spectral_sequence} again. We obtain the theorem by inspection of the spectral sequence.

  For finite fields, by Theorem \ref{th:vdk_gln}, for $2j\leq n+1-p$, we have $L^2_{i,j}=0$. So $p+q\leq k+1$ and $q\leq k$ imply $p+2q\leq 2k+1\leq n+1$, so $E^1_{p,q}=0$, and we obtain the result.
\end{Proof}
\bibliographystyle{alpha}
\bibliography{../biblio}
\end{document}